\documentclass[12pt]{article}
\usepackage{amssymb,latexsym}
\usepackage{amsmath}
\usepackage[dvipdfm]{graphicx}
\usepackage{tikz}
\usepackage{verbatim}
\usetikzlibrary{hobby,backgrounds,calc,trees}
\usepackage{verbatim}
\usetikzlibrary{calc,shapes}
\usetikzlibrary{decorations.pathreplacing,angles,quotes}
\usepackage{rawfonts}
\usepackage{amsmath, psfrag}
\usepackage{amsfonts}
\usepackage{amsthm}
\usepackage[pagebackref]{hyperref}

\newtheorem{theorem}{Theorem}
\newtheorem{observation}[theorem]{Observation}
\newtheorem{lemma}[theorem]{Lemma}
\newtheorem{corollary}[theorem]{Corollary}
\newtheorem{definition}{Definition}

\newtheorem{proposition}[theorem]{Proposition}

\def\qed{\ifhmode\unskip\nobreak\fi\quad\ifmmode\Box\else$\Box$\fi}

\def\qed{\hfill $\Box$}

\renewcommand*{\backref}[1]{}
\renewcommand*{\backrefalt}[4]{%
    \ifcase #1
    Not cited%
    \or
    $\uparrow$#2%
    \else
    $\uparrow$#2%
    \fi
}

\title{Partitioning perfect graphs into comparability graphs}

\author{
{\sl Andr\'as Gy\'arf\'as}\thanks{Supported in part by NKFIH Grant No. K-132696.}\\
\small Alfr\'ed R\'enyi Institute of Mathematics, Budapest, Hungary \\
\texttt{gyarfas@renyi.hu}
\and
{\sl M\'arton Marits}\\
\small Budapest University of Technology and Economics, Budapest, Hungary and \\
\small Yokohama National University, Yokohama, Japan \\
\texttt{marits.marton@gmail.com}
\and
{\sl G\'eza T\'oth}\thanks{Supported by National Research, Development and Innovation Office, NKFIH,
K-131529, ADVANCED 152590, and ERC Advanced Grant "GeoScape" 882971.}\\
\small Alfr\'ed R\'enyi Institute of Mathematics, Budapest, Hungary and\\
\small Department of Computer Science and Information Theory\\
\small Budapest University of Technology and Economics, Budapest, Hungary\\
\texttt{geza@renyi.hu}}

\begin{document}
\maketitle

\begin{abstract} We study how many comparability subgraphs are needed to partition
  the edge set of a perfect graph. We show that many classes of perfect graphs can be partitioned
  into (at most) {\em two} comparability subgraphs and this holds for almost all perfect graphs.
  On the other hand, we prove that for interval graphs
  an arbitrarily large number of  comparability subgraphs might be necessary.
  In particular, for the interval graph defined by the  ${n\choose 2}$ intervals of $[n]$ with integer endpoints,
  we need
  $\Omega(\log\log(n))$ and $O((\log\log(n))^2)$ comparability graphs for the partition.
  \end{abstract}

\section{Introduction}

Perfect graphs are the graphs $G$ in which all induced subgraphs $H\subseteq G$
satisfy the property $\omega(H)=\chi(H)$ where $\omega(H)$ and $\chi(H)$ denote
the size of the largest complete subgraph and the chromatic number, respectively.
Berge put forward two conjectures, the weaker one was that complements of
perfect graphs are also perfect.
This means that the induced subgraphs $H$ also satisfy $\alpha(H)=\theta(H)$,
where $\alpha(H)$ is the cardinality of the largest independent set and
$\theta(H)$ is the minimum number of cliques covering $V(H)$.
The conjecture was proved by Lov\'asz \cite{LO} (weak perfect graph theorem).

The stronger conjecture was that perfect graphs can be characterized
by excluding odd cycles of length at least $5$
and their complements as induced subgraphs.
These graphs are called later {\em Berge graphs}.
This was proved by Chudnovsky, Robertson, Seymour and Thomas \cite{CRST}
(strong perfect graph theorem).
For a survey of many interesting classes of perfect graphs, see \cite{H}.

Here we focus on a well known subclass of perfect graphs, the {\em comparability graphs},
which can be obtained from a partially ordered set as follows.
Vertices correspond to the elements of the partially ordered set,
and two vertices are connected by an undirected edge if and only if the corresponding elements
are comparable. The perfectness of comparability graphs are established in \cite{D50}, \cite{M71}.

The problem we explore is the following:
{\em how many comparability subgraphs are needed to partition or cover the edge set of a perfect graph $G$?}
See the next paragraph for the precise formulation.
The question is very interesting for general graphs instead of perfect graphs,
but here we concentrate on perfect graphs.

For any graph $G$, $E(G)$ denotes its edge set.
Let $p(G)$, denote the minimum number $m$ such that
the edge set of $G$ can be {\em partitioned} into $m$ comparability graphs, that is,
$E(G)=\cup_{i=1}^m E(G_i)$ where the $G_i$-s are edge-disjoint comparability subgraphs of $G$.

Dropping  the condition that the $G_i$-s are edge-disjoint, we get the following definition.
Let $c(G)$, denote the minimum number $m$ such that
the edge set of $G$ can be {\em covered} by $m$ comparability graphs, that is,
$E(G)=\cup_{i=1}^m E(G_i)$ where the $G_i$-s are comparability subgraphs of $G$.

Clearly, $p(G)\ge c(G)$ for every graph $G$,
but it seems not easy to find examples where the gap between them is large.
In fact, we do not know any perfect $G$ for which $p(G)>c(G)$.

The parameter $c(G)$ relates to $\chi$-bounded families of graphs (\cite{GY}, chapter 5).

Covering graphs by comparability graphs is a fundamental tool in combinatorial geometry.
Larman, Matou\v sek, Pach and T\"or\H ocsik \cite{LMPT} observed, that the disjointness graph
of convex sets in the plane can be covered by {\em four} comparability graphs.
As a consequence, they proved that any set of $n$ planar convex sets contain
$n^{1/5}$ pairwise disjoint of pairwise crossing, which is much better
than the bound we would get from Ramsey's theorem.
This general method and observation,
has been used in many contexts \cite{CPS16}, \cite{PRT21}, \cite{RW19}, \cite{PTT24}.


There are triangle-free graphs $G$ with arbitrarily large chromatic number (for a survey see \cite{SS20}).
Since we have $\chi(G)\le \omega(G)^{c(G)}$,
$c(G)$ and $p(G)$ are also arbitrarily large.
However, this argument fails for perfect graphs,
therefore,
it is reasonable to ask whether $c(G)$ and $p(G)$ are  bounded
for perfect graphs.

Throughout this paper
$\log$ denotes the base $2$ logarithm.
An old result of Harary, Hsu and Miller \cite{HHM} (rediscovered many times)  is the following.

\begin{theorem}[\cite{HHM}]\label{HHM}
  The minimum number of bipartite graphs needed to partition the edge set of a graph $G$ is
  $\lceil\log(\chi(G))\rceil$.
\end{theorem}

This result in our context (studying $p(G)$ and $c(G)$ for perfect graphs), easily implies the following.
\begin{corollary}\label{genupbound}  Every perfect graph $G$ can be partitioned into
\begin{itemize}
\item (i) $\lceil\log(\omega(G))\rceil$ bipartite graphs
\item (ii) $1+\lceil\log(\alpha(G))\rceil$ comparability graphs
\end{itemize}
\end{corollary}

\noindent{\bf Proof. }
Since  $\chi(G)=\omega(G)$
for any perfect graph $G$, Theorem \ref{HHM} implies that we can partition $E(G)$ into
$\lceil\log(\omega(G))\rceil$ bipartite graphs, proving (i).
By the weak perfect graph theorem, $\bar{G}$ is also perfect, so
$\chi(\bar{G})=\omega(\bar{G})$. That is, $V(G)$ can be partitioned
into $\alpha(G)$ parts such that each part induces a complete subgraph. The union of these complete subgraphs
is {\em one}
comparability graph. Deleting the edges of these complete subgraphs,
we get a graph $H$ with $\chi(H)\le \alpha(G)$. Applying Theorem \ref{HHM}
for $H$, we get at most $\lceil\log(\alpha(G))\rceil$ further
(bipartite) comparability graphs, proving (ii). \qed

\section{Results}

\subsection{The asymptotic of $p(G)$ and $c(G)$}

Trivially, $p(G)=c(G)=1$  if and only if $G$ is a comparability graph.
It is well-known that there are many perfect graphs that are not comparability graphs.
For those, we have $p(G), c(G)\ge 2$. However, it is not obvious whether there are perfect graphs $G$ with $p(G), c(G)>2$.
In fact, we show the following.

\begin{theorem} \label{twocover} We have $c(G)\le p(G)\le 2$ for almost all perfect graphs.
\end{theorem}

In Theorem \ref{twocover} we use the term ``almost all'' in the following sense.
Let $PE(n)$ denote the set of labelled perfect graphs on $n$ vertices and let $P2PE(n)$ denote the
set of labelled perfect graphs on $n$ vertices that can be partitioned into (at most) two comparability graphs.
Then  Theorem \ref{twocover} states that
$$\lim_{n\rightarrow \infty} {|P2PE(n)|\over |PE(n)|}=1.$$

Theorem \ref{twocover} might suggest that {\em all} perfect graphs can be covered
by, or partitioned into
a constant number (perhaps two) comparability graphs. The proof of Theorem \ref{twocover}
is based on two deep results, the strong perfect graph theorem \cite{CRST} and Theorem \ref{promelsteger}
of Pr\"omel and Steger \cite{PS} below on {\em generalized split graphs}, shortly GSP graphs, defined as follows.

\begin{definition} A graph $G$ is a GSP if either $G$ or its complement can
  be written in the following form:
  the vertex set is partitioned into $V_1,V_2$ where $V_1$
  is a complete graph and $V_2$ is the union of vertex disjoint
  complete graphs with no edges between them.
\end{definition}

It is important for the proof of Theorem \ref{twocover}
that GSP graphs are perfect. When GSP graphs were defined in
\cite{PS}, the proof of their perfectness was based
on the incorrect claim that they form a subclass of weakly
triangulated graphs (the six-cycle $C_6$
is a not weakly triangulated GSP graph).
Later in \cite{EW11} the perfectness of GSP graphs was derived from
the strong perfect graph theorem \cite{CRST}.
In fact, the weak perfect graph theorem of Lov\'asz \cite{LO}
is sufficient: it reduces the proof to show that
$\alpha(G)=\theta(G)$ for graphs $G$ in the form
$V(G)=V(K_0)\cup V(K_1)\cup \dots \cup V(K_t)$, where $K_i$-s are complete
graphs and there are no edges between
$K_i,K_j$ if $0<i<j\le t$. 
The statement is obvious when $\alpha(G)=t+1$.
If $\alpha(G)=t$ then every vertex of $K_0$ is adjacent to all vertices
of $K_i$ for some $1\le i\le t$, implying that the $t$-cover
$\cup_{i=1}^t V(K_i)$  can be extended by $V(K_0)$ to a $t$-cover of $V(G)$.

\begin{theorem}\label{promelsteger}[\cite{PS}] Almost all Berge graphs are GSP.
\end{theorem}
As in Theorem \ref{twocover}, ``almost all'' in Theorem \ref{promelsteger} means
$$\lim_{n\rightarrow \infty} { |GSP(n)|\over|BE(n)|}=1,$$
where $BE(n),GSP(n)$ denote the set of all labelled Berge graphs (respectively GSP graphs) on $n$ vertices.
Now we know, by the strong perfect graph theorem \cite{CRST},
that Berge graphs are exactly the perfect graphs, so $BE(n)=PE(n)$.

The link we need towards the proof of Theorem \ref{twocover} is the following property of GSP graphs.

\begin{theorem} \label{gensplit} For any GSP graph $G$, $p(G)\le 2$.
\end{theorem}

\noindent
    {\bf Proof. } Assume first that the vertex set of a GSP graph $G$ is partitioned into $V_1,V_2$ where $V_1$
    is a complete graph and $V_2$ is the union of vertex disjoint complete graphs $K_i$ (with no edges between them).
    Then the complete subgraphs spanned by $V_1$ and by the $V(K_i)$-s form a graph of vertex disjoint cliques,
    which is a comparability graph
    (just take a linear order on the vertex set of each clique, with  pairs from different cliques not comparable),
    and the other comparability graph is the bipartite graph
    $G[V_1,V_2]$.

    In the second (complementary) case $V_1$ is an independent set in $G$ and $V_2$ spans a
    complete partite graph $H$, which is obviously a comparability graph.
    Then $H$ and the bipartite graph spanned by $[V_1,V(H)]$ in $G$ give the required partition.   \qed

\medskip

\noindent
{\bf Proof of Theorem \ref{twocover}.} Using Theorem  \ref{gensplit} and the strong perfect graph theorem, we have
$GSP(n)\subset P2PE(n)$ and $BE(n)=PE(n)$, therefore from Theorem \ref{promelsteger},

$$1=\lim_{n\rightarrow \infty} {|GSP(n)|\over |BE(n)|}\le \lim_{n\rightarrow \infty}{|P2PE(n)|\over |PE(n)|}\le 1,$$
implying that the second limit is also equal to $1$, finishing the proof. \qed

\subsection{Perfect graph classes with $p(G)\le 2$}

The {\em line graph} of a graph $G$ is obtained by representing the edges of $G$ by vertices and
considering two such vertices adjacent if the corresponding edges of $G$ have a common vertex.
The line graphs of bipartite graphs, LBIP, form a well-known class of perfect graphs \cite{G04}.

\begin{theorem}\label{bipline} For $G\in $LBIP, $p(G)\le 2$.
\end{theorem}

\medskip

\noindent{\bf Proof. }  Let $G$ be the line graph of a bipartite graph $H=[A,B]$.
The edges of $G$ are labelled with $a$ ($b$) if the the corresponding edges of $H$ intersect in $A$ (in $B$).
Clearly, the edges labelled with $a$ (with $b$) are unions of complete graphs with no edges between them,
thus they give the required partition of $E(G)$ into two comparability graphs. \qed

\medskip
Some perfect graphs are defined as intersection graphs of geometric objects.
A graph $G$ is an {\em interval graph} if it is the intersection graph of a system of closed intervals of a line.
If the intervals have length one, their intersection graph is a {\em unit interval graph}.
Interval graphs are perfect \cite{G04},
have several characterizations, the most relevant in our case is from Gilmore and Hoffman \cite{GH}:
a graph is an interval graph if and only if it does not contain an induced $C_4$ and its complement is a comparability graph.
A more general family of perfect graphs is the intersection graphs of subtrees of a tree.
It is characterized (independently in \cite{GAV}, \cite{WA}) as the family of graphs that do
not contain induced cycles of length at least four and called {\em chordal, alias triangulated graphs}.


\begin{theorem}\label{unit} For unit interval graphs $G$, $p(G)\le 2$.
\end{theorem}

\noindent{\bf Proof. }  Let $U$ be a system of unit intervals and partition them as follows.
Let $I_1=[a_1,b_1] \in U$ be the interval with leftmost right endpoint and define $X_1$ as the set of intervals of
$U$ containing $b_1$.
If $I_1,\dots, I_k$ and $X_1,\dots,X_k$ are defined and there are further intervals, select from
them the interval $I_{k+1}=[a_{k+1},b_{k+1}]$
with the leftmost right endpoint.
Define $X_{k+1}$ as the set of intervals not in $\cup_{i=1}^k X_i$ and containing $b_{k+1}$.
The procedure ends with $[a_m,b_m]$ and $X_m$ when all intervals of $U$ are partitioned.

The unit interval graph of $U$, $G_U$, can be partitioned into two comparability graphs as follows.
Define $H_1$ as the edge set of $G_U$ defined by two intervals in  the {\em same}
$X_i$, for any $i$.
The graph $G_1(U, H_1)$ is the union of $m$ vertex-disjoint complete graphs.
Let $H_2$ be the set of remaining edges, i.e. those that correspond to
intersecting from {\em different}
$X_i$-s.
Since $U$ contains unit intervals,
edges of the graph $G_2(U, H_2)$ must represent intersecting intervals belonging to
$X_i$ and $X_{i+1}$ for some $i$.
Thus $G_2(U, H_2)$ is a bipartite graph.
Clearly, $G_1$ and $G_2$ are both comparability graphs, finishing the proof.   \qed

\medskip

\begin{figure}[h]
	\begin{center}
		\includegraphics[scale=0.45]{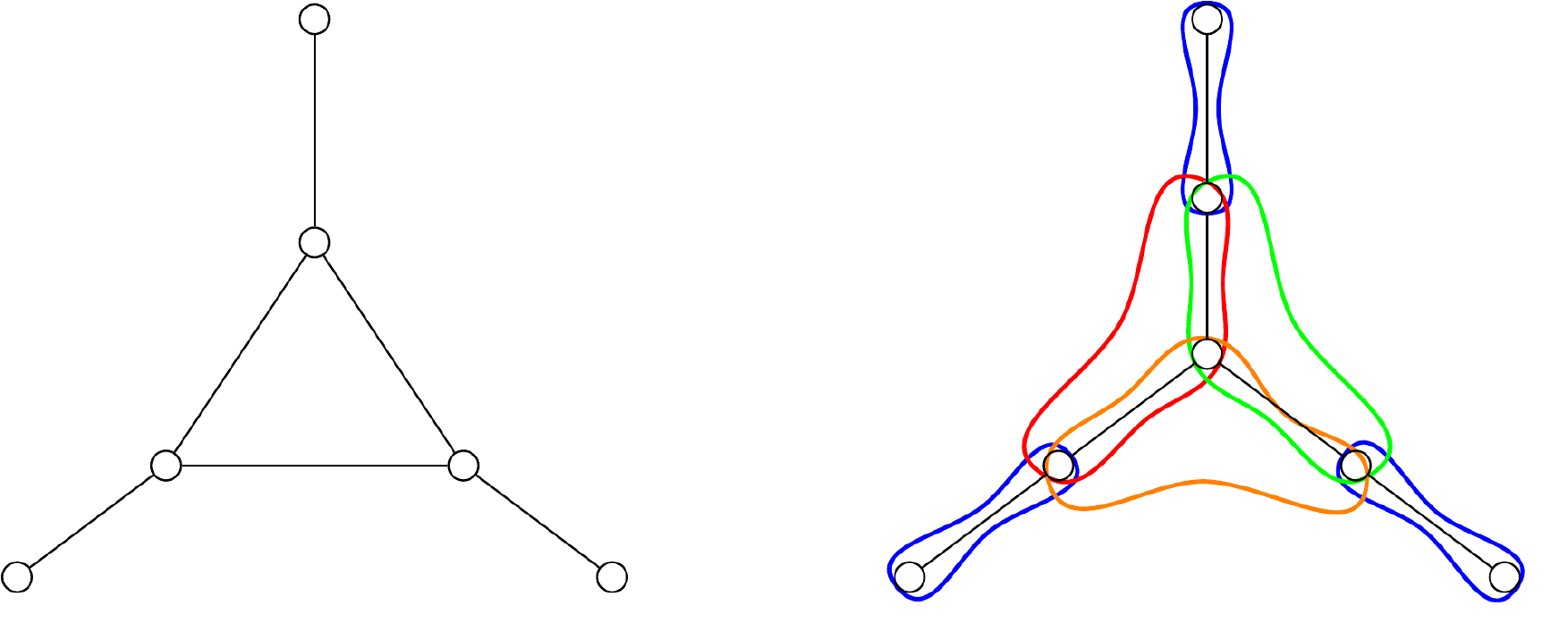}
	\end{center}
	\caption{A graph that is not a comparability graph, represented as the disjointness graph of subtrees of a tree} \label{fig.noncomp}
\end{figure}

The complements of interval graphs, i.e. the disjointness graphs of interval systems are comparability graphs.
Consider a more general class, disjointness graphs of subtrees.
Let $T$ be a tree. The vertices of the graph $G$ correspond to some subtrees of $T$,
two vertices are connected if the corresponding subtrees are disjoint.
These graphs are also called {\em co-triangulated}.
These graphs are not always comparability graphs, Figure \ref{fig.noncomp} shows an example.

Nevertheless, they can be partitioned into two comparability subgraphs.
\begin{theorem} \label{co-triangulated}
  Assume that $G$ is a co-triangulated graph,
  i.e. the disjoint\-ness graph of subtrees of a tree. Then $c(G)\le p(G)\le 2$.
\end{theorem}

\medskip
\noindent{\bf Proof. } Let $T$ be a tree of $n$ vertices and $\cal{F}$ is a system of subtrees of $T$.
We construct two comparability graphs on the disjoint pairs of $\cal{F}$.
 Let $r$ be any vertex of $T$, the root. For any subtree $A\in {\cal{F}}$
there is a unique vertex of $A$ that is closest to $r$, call this vertex the root
of $A$.
For the first partial order $\prec_1$,
let $A,B\in {\cal{F}}$
be two disjoint subtrees. We say that
$A\prec_1 B$ if the unique path from the root of $B$ to $r$ contains the root of $A$. Intuitively,
we can say that $A\prec_1 B$ ``if $A$ is closer to $r$ than $B$''.
It is clear that $\prec_1$ is antisymmetric. Suppose that  $A\prec_1 B$ and  $B\prec_1 C$ for subtrees $A$, $B$, $C$.
Then the path $P_C$ from the root of $C$ to $r$ contains the root of $B$. But then $P_C$ contains the path $P_B$
from the root of $B$ to $r$. Since $P_B$ contains the root of $A$, $P_C$ also contains
the root of $A$. Moreover, the root of $B$ separates $A$ and $C$, that is, if we remove the root of $B$ from $T$,
$A$ and $C$ are in different components,
so they are disjoint.
Therefore, $A\prec_1 C$, consequently, $\prec_1$ is transitive, so it is
a partial ordering on the disjoint pairs.

For any $i>0$, let  $N^{i}$ be the set of $i$-tuples of natural numbers, and let
Let $N^{\le n}=\cup_{i=0}^nN^i$.
We define a partial order $\prec_{lex}$ on $N^{\le n}$ as follows.

Suppose that  $\alpha\in N^{i}$  $\beta\in N^{j}$.
We say that  $\alpha$ is a prefix of  $\beta$, if
$i\le j$
and we can get  $\alpha$ from  $\beta$ by deleting its last $j-i$ digits.

For $\alpha, \beta\in N^{\le n}$, let $\alpha\prec_{lex}\beta$
if none of $\alpha$ and $\beta$ is a prefix of the other,  and
$\alpha\prec\beta$ in the lexicographic ordering.

For the second partial order on the disjoint trees,
we first define the labeling function $\psi : V(T) \rightarrow N^{\le n}$
in the following recursive way. In the root, we define $\psi= \emptyset$. For neighbors
of the root, we arbitrarily choose an ordering, and assign the numbers
$\{0,\ldots, k\}$
where $k$ is the degree of $r$ minus one.
In general, if a vertex $v$ has the label $\psi(v)$
and $v$ has at least one child, they get labels $\psi(v)0, \dots \psi(v)k$ where
$k$ is the number of children minus one.

The value of $\psi$ at each
vertex is a sequence of natural numbers.
Observe that for disjoint subtrees $A,B$ with roots $a,b$, $A \prec_1 B$ if and only if
$\psi(a)$ is a prefix of $\psi(b)$.

Now let $A,B$ be two disjoint subtrees of $T$, with roots $a,b$ respectively,
which are not comparable by $\prec_1$.
Then let $A \prec_2 B$ if and only if  $\psi(a)\prec_{lex} \psi(b)$.

It is easily seen that $\prec_2$ is irreflexive and antisymmetric,
it remains to prove
transitivity. Let $A,B,C$ be subtrees of $T$ with roots
$a,b,c$. If $A\prec_2 B$ and
$B \prec_2 C$, then clearly, $\psi(a) \prec_{lex}\psi(b)\prec_{lex}\psi(c)$.
This  implies that none of $\psi(a)$ and $\psi(c)$ is a prefix of the other and
$\psi(a) \prec_{lex} \psi(c)$.

Moreover, for any vertex $u$ of $A$, $\psi(a)$ is a prefix of $\psi(u)$,
therefore, $\psi(u)\prec_{lex}\psi(b)$. On the other hand,
for any vertex $v$ of $C$, $\psi(c)$
is a prefix of $\psi(v)$, therefore, $\psi(b)\prec_{lex}\psi(v)$.
Consequently,
$A$ and $C$
must be disjoint.

Therefore, the comparability graphs of $\prec_1$ and $\prec_2$ on ${\cal F}$
form
a partition of the disjointness graph on  ${\cal F}$.   \qed

\subsection{Interval graphs}

In spite of what the previous results suggest, there are perfect graphs
$G$ with arbitrary large $c(G)$, consequently, arbitrarily large $p(G)$ as well.
This happens for interval graphs, a well-known family of perfect graphs defined
as intersection graphs of intervals of a line (see Theorem \ref{loglog} below).
For $n\ge 2$, let $G_n$ be the interval graph defined as the intersection graph
of the $n\choose 2$ closed intervals with endpoints among $n$ consecutive positive integers,
called {\em base points}. Note that $|V(G_n)|={n\choose 2}$.
For convenience, set $$c(n)=c(G_n),p(n)=p(G_n).$$
Since $G_2$ has one vertex and no edge, we set $c(2)=p(2)=0$.

Two distinct closed intersecting intervals $I_1,I_2$ with
integer endpoints are called {\em nested } if one contains the other.
If $I_1,I_2$ are non-nested then they are either {\em crossing}
(if  $|I_1\cap I_2|\ne 1$) or {\em touching} (if  $|I_1\cap I_2|=1$).

An edge of $G_n$ is called {\em nested}, {\em crossing}, or {\em touching}
if the corresponding
pair of
intervals is {nested},
{crossing}, or {touching}.

Let $S_n$ denote the subgraph of touching edges of $G_n$.
The graph $S_n$ is often referred as the {\em shift graph}.
It is well-known that $\chi(S_n)=\lceil\log(n)\rceil$
(in \cite{TR} this is attributed to Andr\'as Hajnal). From Theorem  A. we have the following.

\begin{corollary}\label{shift}
  $p(S_n)\le \lceil \log(\chi(S_n))\rceil = \lceil\log(\lceil\log(n)\rceil)\rceil$
  and the partition can be done by bipartite graphs.
\end{corollary}

The {\em double shift graph}, $DS_n$, is a close relative of shift graphs,
introduced by Erd\H os and Hajnal \cite{EH}.
The $n\choose 3$ vertices of $DS_n$ correspond to triples $(i, j, k)$, $1\le i<j<k\le n$.
Two vertices, $(i, j, k)$ and $(x, y, z)$ are connected by an edge if and only if
(i) $j=x$ and $k=y$ or (ii) $i=y$ and $j=z$. From a result of Erd\H os and Hajnal
\cite{EH} it follows that
\begin{corollary}\label{doblesh} $\chi(DS_n)\ge\log\log(n)$.
\end{corollary}

\begin{theorem}\label{loglog} $c(n)\ge  {1\over 2}\log\log(n)$.
\end{theorem}

\medskip
\noindent{\bf Proof. }  Assume that $G_n$ is the union of $t$ comparability subgraphs,
$H_1,\dots,H_t$. Consider each $H_i$ with a transitive orientation.
     Assign a type to each touching edge $e=\{[i,j],[j,k]\}$
    as a pair $(m,+)$ or $(m,-)$ where the index $m$ indicates that the touching edge is
    in $H_m$ and the + or - indicates whether the touching edge is oriented from $[i,j]$
    to $[j,k]$ or backwards.
    (Since $H_1,\dots,H_t$ is a cover, we might have several choices for $m$, in this case
    take an arbitrary one.) The number of types is clearly $2t$.

    Note that the types of the touching edges $e=\{[i,j],[j,k]\}$ and $f=\{[j,k],[k,l]\}$ must be different
    for all $1\le i<j<k<l\le n$, otherwise the transitivity assumption implies that there is an edge of $G_n$
    from $[i,j]$ to $[k,l]$ - a contradiction since those intervals are disjoint.

    Therefore, we obtained a proper coloring of the double shift graph $DS_n$ by $2t$ colors:
    the color of the vertex $(i, j, k)$ corresponds to the type of
    the touching edge $e=\{[i,j],[j,k]\}$. Consequently, $2t\ge \chi(DS_n)\ge \log\log(n)$, finishing the proof.  \qed

\medskip

A natural upper bound of $p(n)$ comes from Theorem \ref{HHM}.
Indeed, $\chi(G_n)=\omega(G_n)\sim {n^2\over 4}$ because
any set of pairwise crossing intervals have a comon point, therefore,
$\omega(G_n)$ is the
maximum number of intervals containing one of the $n$ base points.
Then Theorem \ref{HHM} gives a $\sim 2\log(n)$ upper bound for $p(n)$
(partitioning into bipartite graphs).
Here we prove an $O((\log\log(n))^2)$  upper bound, closer to the lower bound of Theorem \ref{loglog}.

\begin{theorem}\label{ub} For any $n\ge 4$ we have
  $$p(n)\le {1\over 2}(\lceil\log(\lceil\log(n)\rceil)\rceil)^2+{5\over 2}\lceil\log(\lceil\log(n)\rceil)\rceil-1.$$
\end{theorem}

Before the proof, we need some preparations.
Let $F_n$ be the graph obtained from $G_n$ by removing all touching edges of $G_n$,
i.e. edges of $S_n$ (without removing their endpoints).
Let
$p'(n)=p(F_n)$.

\begin{lemma}\label{rec}  We have $p'(4)=1$. For $n\ge 3$ we have
 $$p'(n^2)\le  p'(n)+ \lceil\log(\lceil\log(n)\rceil)\rceil + 2.$$
\end{lemma}

\begin{proof}[Proof of Lemma \ref{rec}]
It is easy to check that $p'(4)=1$, i.e. $F_4$ is a comparability graph.
Indeed, in the vertex ordering
$$(2,4),(1,2),(3,4),(2,3),(1,3),(1,4),$$

orienting all the ten edges of $F_4$ forward, we get a transitive orientation, thus $F_4$ is indeed a comparability graph.
Let $n\ge 3$ now.
Assume that
$P=H_1,\dots,H_{p'(n)}$ is a partition of the edges of $F_n$ into $p'(n)$
comparability subgraphs. Partition the base point set $[n^2]$ of $F_{n^2}$ into $n$ consecutive blocks of size
$n$,  $B_1,\dots,B_n$, where for $i\in [n]$,

$$B_i=\{(i-1)n+1,\dots,in\}.$$

The intervals (that is, vertices of
$F_{n^2}$)
with both endpoints in the same block $B_i$
are called {\em inner} intervals.
The intervals with endpoints in different blocks $B_i,B_j,
i\ne j$  are the {\em outer} intervals.

Let $e$ be the edge of $F_{n^2}$ defined by a pair of
outer intervals $[a, b]$ and $[c, d]$,
where $a\in B_{a'}$,  $b\in B_{b'}$,   $c\in B_{c'}$,   $d\in B_{d'}$.
Then $e$ is called {\em outer-nested, outer-crossing, outer-touching, or outer-degenerate}
if the intervals $[a',b'],[c',d']$ are nested, crossing, touching or coincide (i.e. $a'=c',b'=d'$).

For $i\in [n]$, let $F(B_i)$ denote the subgraph of $F_{n^2}$ induced by the inner intervals of $B_i$.

\begin{definition}
  The {\em blow-up of type (a)} of a graph $G$ is the following.
  Replace vertices  $v_i$ of $G$ by pairwise disjoint vertex sets $V_i$.
  For each edge $e=(v_i,v_j)$ of $G$, add a complete bipartite graph between
  the corresponding sets $V_i$ and $V_j$.

  The {\em blow-up of type (b)} of a graph $G$ is the following.
Replace vertices  $v_i$ of $G$ by pairwise disjoint vertex sets $V_i$.
  For each edge $e=(v_i,v_j)$ of $G$, add an arbitrary set of edges between
  the corresponding sets $V_i$ and $V_j$.
\end{definition}

The following observation is obvious.

\begin{observation}\label{blow}
  The blow-up of type (a) of a comparability graph is a comparability graph.
    The blow-up of type (b) of a bipartite graph is a bipartite graph.
\end{observation}

Returning to the proof of Lemma \ref{rec}, $F_{n^2}$ is partitioned as follows.

\begin{itemize}

\item (i)  Clearly, $F(B_i)$ is isomorphic to $F_n$ so it has a partition
  $P_i= H^i_1, \ldots, H^i_{p'(n)}$ into $p'(n)$ comparability graphs. Since unions of vertex
  disjoint comparability graphs are comparability graphs,
  $P^*= H^*_1, \ldots, H^*_{p'(n)}$
is a partition on the edges of $F_{n^2}$ defined by
inner intervals into comparability graphs $H^*_j=\cup_{i=1}^nH^i_j$.

Let $NC\subseteq F_{n^2}$ be a graph whose vertices correspond to the outer intervals and edges
  correspond to outer-nested and
  outer-crossing pairs of intervals.
  Clearly, the graph $NC$ is a blow-up of type (a) of $F_n$. Therefore, by Observation \ref{blow}, there is a partition
  of its edges into $p'(n)$ comparability graphs,
 $P'=H'_1,\dots,H'_{p'(n)}$.
Again, since  $P'$ is defined on
outer intervals while $P^*$ is on inner intervals and
   unions of vertex
   disjoint comparability graphs are comparability graphs, we can take the union of the partitions
   $P^*$ and $P'$ to get
   a partition of the {\em edges defined by inner intervals and by outer-nested and outer-crossing pairs} into
   $p'(n)$ comparability graphs.

 \item (ii) Consider now the graph $OS$ of {\em intersecting outer-touching intervals} of $F_{n^2}$.
   An outer-touching edge is a pair of outer-touching intervals $[i,j], [k,l]$ $i<j$,$k<l$ where
   $i\in B_a$, $j, k\in B_b$ and $a<b<c$, $l\in B_c$, so it corresponds to a touching pair $[a, b], [b, c]$. Therefore,
   $OS$ is a blow-up of type (b) of $S_n$. Observe, that a pair of outer-touching intervals are not necessarily intersecting,
   that's why we have a blow-up of type (b) and not (a).

   Corollary \ref{shift} gives a partition of the edges of $S_n$ defined by touching pairs of intervals
   into $\lceil\log(\lceil\log(n)\rceil)\rceil$ bipartite graphs.
   By Observation \ref{blow}, it corresponds to a partition of the edges of $OS$
   into the same number of bipartite graphs. Since bipartite graphs are comparability graphs, we
   obtained a partition of $OS$ into  $\lceil\log(\lceil\log(n)\rceil)\rceil$ comparability graphs.

 \item (iii) Let $OD$ be the graph of {\em intersecting outer-degenerate intervals} of $F_{n^2}$.
   For any interval $[i, j]$, $1\le i<j\le n^2$, $i\in B_a$, $j\in B_b$, the type of $[i,j]$ is the pair $(a,b)$.
   Observe that two outer intervals form an intersecting outer-degenerate pair if and only if they are of the same type.
Therefore, the graph $OD$ is the disjoint union of ${n\choose 2}$ cliques, which is a comparability graph.

\item (iv) Finally, there are the edges of $F_{n^2}$ that correspond to a crossing or nested pair of an inner and an outer interval.
But these edges form a bipartite graph, so it is another comparability graph.

\end{itemize}

Since all edges of $F_{n^2}$ are partitioned in (i)-(iv), summing up, we get the statement of Lemma \ref{rec}.
\end{proof}

A repeated application of Lemma \ref{rec} gives the following.

\begin{corollary}\label{p'}
  For any  $n= 2^{2^k}$ with $k\ge 1$ we have
  $$p'(n)=p'(2^{2^{k}})\le p'(2^{2^{k-1}})+k-1+2 \le \dots \le {k\choose 2}+2k-1.$$
\end{corollary}

Now we are ready to prove Theorem \ref{ub}.

\begin{proof}[Proof of Theorem \ref{ub}]
It is enough to prove the statement for $n=2^{2^k}$ where $k\ge 2$.
Since $p(n)$ is non-decreasing, this implies the statement for $2^{2^{k-1}}<n < 2^{2^k}$ as well.

Let $n=2^{2^k}$. Observe that the edge set of $G_n$ is the disiont union of the edge set of $S_n$ and $F_n$.
Combining the partitions given by Corollary \ref{shift} and \ref{p'} respectively we get that

$$p(n) \le p'(n)+ p(S_n)\le {k\choose 2}+2k-1+ \lceil\log(\lceil\log(n)\rceil)\rceil={k\choose 2}+3k-1$$
$$\le {1\over 2}(\lceil\log(\lceil\log(n)\rceil)\rceil)^2+{5\over 2}\lceil\log(\lceil\log(n)\rceil)\rceil-1.$$
finishing the proof of Theorem \ref{ub}.
\end{proof}

\subsection{A ``small'' perfect graph with $c(G)=p(G)=3$}

Theorem \ref{loglog} gives huge perfect graphs, for example the interval graph $G_n$
with $c(G_n)>2$  has more than
$2^{16}\choose 2$ vertices.
For the $c(G)=3$ case we can construct a smaller
example  based on the cartesian product  $K_n\square K_m$, the line-graph of the complete bipartite graph $K_{n,m}$.
Its vertices are $(i, j)$, $1\le i\le n$, $1\le j\le m$,
two different vertices, $(i, j)$ and $(i', j')$ are connected by an edge if and only if
$i=i'$ or $j=j'$.
Let $H_{9,4}$ be the graph with $72$ vertices,
obtained from $K_9\square K_4$ by attaching a pendant edge at each vertex.

\begin{proposition}\label{H9} The graph $H=H_{9,4}$ is perfect and $p(H)=c(H)=3$.
\end{proposition}

\medskip

\noindent{\bf Proof.}

Let $V(H)=K\cup W$ where
$$K=\{v_{i,j}: 1\le i\le 9, 1\le j\le 4\}, W=\{w_{i,j}:  1\le i\le 9, 1\le j\le 4\}.$$
For $1\le j\le 4$, let $X_j=\{v_{i, j} : 1\le i \le 9\}$ be the {\em $j$-th row} of $K$ and
for  $1\le i\le 9$, let $Y_i=\{v_{i, j} : 1\le j \le 4\}$ be the {\em $i$-th column} of $K$.

There are three types of edges,
(i)  for  $1\le i\le 9$, $1\le j\le 4$,  $e_{i,j}=\{v_{i,j},w_{i,j}\}$ are the pendant edges of $H$,
(ii) for $1\le j\le 4$, $X_j$ induces a complete graph
(iii)  for $1\le i\le 9$, $Y_i$ also induces a complete graph.

  The graph $H$ is perfect since $K_9\square K_4$ is the line graph of the complete bipartite graph
  $K_{9,4}$ and adding pendant edges to a perfect graph preserves perfectness.

  To see that $c(H)\le p(H)\le 3$,
  note that the transitive orientations on the complete graphs spanned by the
  $X_i$-s (resp. $Y_i$-s), and any orientation on the pendant edges provide a partition of $E(H)$ into three comparability subgraphs.

  On the other hand, suppose for a contradiction that
  $H$  is the union of two comparability subgraphs, $A$ and $B$.
Consider $A$ and $B$
with their transitive orientations (which might introduce multiply directed edges).

A vertex $v_{i,j}\in X_j$ is  $A$-good (resp. $B$-good)
if it is
the midpoint of an oriented two-edge path in
row $X_j$ with both edges in $A$ (resp. $B$).
A vertex $v_{i,j}$ is good if it is $A$-good or $B$-good and bad otherwise.

Observe that if $v_{i,j}$ is bad, than it has to be maximal or minimal in both
transitively oriented subgraphs $A$ and $B$, restricted to $X_j$.
So there are four types of bad vertices, ($A$-max, $B$-max),
($A$-max, $B$-min),  ($A$-min, $B$-max),  ($A$-min, $B$-min).

Let $u$ and $v$ be two bad vertices
in a row $X_i$. Suppose that $uv$ belongs to $A$.
Then $u$ and $v$ cannot be $A$-max (resp. $A$-min) at the same time.
Similarly, if $uv$ belongs to $B$,
then $u$ and $v$ cannot be $B$-max (resp. $B$-min) at the same time.
Therefore, in each row $X_i$ there is at most one bad vertex of each type, so there are at most $4$ bad vertices.
Consequently, in each row there are at least $5$ good vertices.
So $K$ contains at least $20$ good vertices. It follows that one of the columns contains at least $3$ good vertices.
Assume without loss of generality that it is $Y_1$ and the three good vertices are $v_{1,1}$, $v_{1,2}$, $v_{1,3}$.

Suppose that $v_{1,1}$ is $A$-good and $v_{1,2}$ is $B$-good.
The edge $v_{1,1}v_{1,2}$ belongs to $A$ or to $B$, and in both cases, by transitivity,
we must have an edge of $H$ which is neither vertical, nor horizontal, which is a contradiction.
Therefore,  $v_{1,1}$ and $v_{1,2}$ are of the same type. Analogously, $v_{1,3}$ is also of the same type, say, all of them are $A$-good.
Consider now the edges between  $v_{1,1}$, $v_{1,2}$ and  $v_{1,3}$. If any of them belongs to $A$, then again  by transitivity,
we must have an edge of $H$ which is neither vertical, nor horizontal, which is a contradiction.
So, all three edges belong to $B$. Since $B$ has a transitive orientation,
$v_{1,1}$, $v_{1,2}$ and  $v_{1,3}$ form an oriented path in $B$. Suppose that $v_{1,1}$ is its middle vertex.
Then $v_{1,1}$ is the middle vertex of an $A$-path and $B$-path in $K$.
Finally, no matter whether $e_{1,1}$ belongs to $A$ or to $B$ and no matter which direction is assigned to it,
by the transitivity of $A$ or $B$ we must have another edge from $w_{1,1}$ which is a contradiction.
\qed.

The perfect graph $H_{9,4}$ in Proposition \ref{H9} is not the smallest of its kind.
N. Alm\'asi \cite{A24} proved that $K_5\square K_5$ with a pendant edge at each vertex also works.
However, there is no chance to build similar examples for $p(G)>3$ using higher dimensional products, because $K_2\square K_2\square K_3$ is not perfect, it contains an induced $7$-cycle, see \cite{KP}.
\medskip

\noindent{\bf Remarks.} The main result of this note is that almost all perfect graphs $G$ we have $c(G)\le p(G)\le 2$.
However, there are some perfect graphs with $c$ and $p$ arbitrarily latge.
We could not find a prefect graphs with $p\neq c$.

\noindent {\bf Problem 1.}
Is it true that for every perfect graph $G$, $p(G)=c(G)$?

Our examples with large $p$ and $c$ have double
exponentially many vertices. We believe that there are much smaller examples.

\noindent {\bf Problem 2.}
Let $f(m)$ be the smallest number of vertices of a perfect graph $G$ with
$c(G)\le m$. Is it true that $f(m)=2^{2^{\Omega(m)}}$?

\noindent{\bf Acknowledgement. }
The authors are very grateful for the helpful suggestions and corrections of the referees.
AI was not used in the preparation of this paper.

\end{document}